\input amstex
\magnification =\magstep 1
\documentstyle{amsppt}
\pageheight{9truein}
\pagewidth{6.5truein}
\NoRunningHeads
\baselineskip=16pt

\topmatter

\title On groups in which  every element has a prime power order
and which satisfy some boundedness condition
\endtitle

\author Marcel Herzog*, Patrizia Longobardi** and Mercede Maj**
\endauthor

\affil *School of Mathematical Sciences \\
       Tel-Aviv University \\
       Ramat-Aviv, Tel-Aviv, Israel
{}\\
       **Dipartimento di Matematica \\
       Universit\`a di Salerno\\
       via Giovanni Paolo II, 132, 84084 Fisciano (Salerno), Italy
\endaffil

\thanks This work was supported by the National Group for Algebraic and
Geometric Structures, and their Applications (GNSAGA - INDAM), Italy.
\endthanks

\abstract  In this paper
we shall deal with periodic groups, in which each element has a prime power order.  A group $G$ will be  called a $BCP$-group
if each element of $G$ has a prime power order and  for each $p\in \pi(G)$ there exists a positive integer
$u_p$  such that each $p$-element of $G$ is of order $p^i\leq p^{u_p}$.
A group $G$ will be called a $BSP$-group if each element of $G$ has a  prime power order and for each $p\in \pi(G)$ there exists a positive integer
$v_p$  such that each finite $p$-subgroup of $G$ is of order $p^j\leq p^{v_p}$. Here $\pi(G)$  denotes the set  of all primes dividing the order of some element of $G$.
Our main results are the following four theorems. Theorem 1: Let $G$ be a finitely generated $BCP$-group. Then $G$ has only a finite number of normal
subgroups of finite index. Theorem 4: Let $G$ be a locally graded $BCP$-group. Then $G$ is a locally finite group. Theorem 7: Let $G$ be a  locally graded $BSP$-group. Then $G$ is a finite group.
Theorem 9: Let $G$ be a $BSP$-group satisfying $2\in \pi(G)$. Then $G$ is a locally finite group.
\endabstract
\endtopmatter
\document

\heading I. Introduction \\
\endheading

In this paper we shall deal with periodic groups, in which each  element has a prime power order. 
The set  of all primes dividing the order of some element of $G$ will be denoted by $\pi(G)$.

In the paper [3] of A.L. Delgado and Y.-F. Wu, groups with each element having a prime power order were called  $CP$-groups. Such groups are of course
periodic. We shall investigate $CP$-groups which satisfy some boundedness condition, as defined below.

\definition {Definitions} A group $G$ will be  called a $BCP$-group if each element of $G$ has a prime power order and  for each $p\in \pi(G)$ there exists a positive integer
$u_p$  such that each $p$-element of $G$ is of order $p^i\leq p^{u_p}$.

A group $G$ will be called a $BSP$-group if each element of $G$ has a prime power order and for each $p\in \pi(G)$ there exists a positive integer
$v_p$  such that each finite $p$-subgroup of $G$ is of order $p^j\leq p^{v_p}$.
\enddefinition

Notice that each $BSP$-group is a  $BCP$-group and each $BCP$-group is a $CP$-group. Moreover, the $BCP$-property and the $CP$-property are  inherited by subgroups and quotient groups, and hence by sections. The $BSP$-property is inherited by subgroups.

The investigation of $BSP$-groups is obviously related to the famous problem that W. Burnside raised in 1902: does a finitely generated group of finite exponent have to be finite? (see [2]). In fact, for any positive integers $n, s$ and every prime $p$, the free Burnside group $B(n, p^s)$ on $n$ generators and of exponent $p^s$ is a $BCP$-group. The knowledge of this problem is very incomplete, for example it is still open if $B(2,5)$ or $B(2,8)$ is finite (see for example [7]).  On the other hand it is well-known that $B(n,e)$ is infinite for sufficiently large exponent $e$ (see [1], [4], [5]). Moreover, A.Yu. Ol'sanskii constructed for any sufficiently large prime $p$ (one can take $p > 10^{75}$) a finitely generated infinite simple group of exponent $p$. (see [8]).

Our aim in this paper  is to find properties of $BCP$-groups and $BSP$-groups, which force these groups to be either finite or locally finite.
Our main results are the following four theorems. Recall that a group $G$ is locally graded if each non-trivial finitely generated subgroup of $G$ has a proper normal subgroup of finite index.

\proclaim {Theorem 1} Let $G$ be a finitely generated $BCP$-group. Then $G$ has only a finite number of normal
subgroups of finite index.
\endproclaim

\proclaim {Theorem 4} Let $G$ be a locally graded $BCP$-group. Then $G$ is a locally finite group.
\endproclaim

\proclaim {Theorem 6} Let $G$ be a  locally finite $BSP$-group. Then $G$ is a finite group.
\endproclaim

\proclaim {Theorem 7} Let $G$ be a  locally graded $BSP$-group. Then $G$ is a finite group.
\endproclaim

We are grateful to the referee of  this paper, for  suggesting  that we consider also  $BCP$-groups and $BSP$-groups $G$, which satisfy
the condition $2\in \pi(G)$. In this direction, we proved the following three additional theorems.

\proclaim {Theorem 5} Let $G$ be a $BCP$-group satisfying $2,3\in \pi(G)$ and suppose that $u_2=1$ and $u_3\in \{1,2\}$. Then
$G$ is a locally finite group.
\endproclaim

\proclaim {Theorem 8} Let $G$ be a $BSP$-$2$-group. Then $G$ is a  finite group.
\endproclaim

\proclaim {Theorem 9} Let $G$ be a $BSP$-group satisfying $2\in \pi(G)$. Then $G$ is a locally finite group.
\endproclaim

The next two sections will deal with $BCP$-groups and $BSP$-groups, respectively.

\heading II. $BCP$-groups \\
\endheading
This sections deals with $BCP$-groups. First we present our  basic result concerning $BCP$-groups.
It is well known  that  finitely generated groups  have only a
finite number of subgroups of a {\it given} finite index.
In particular, each such group has only a finite number of normal subgroups of a {\it given} finite index. We shall  show that
finitely generated $BCP$-groups have only a finite number of normal subgroups of an {\it arbitrary} finite index.

\proclaim {Theorem 1} Let $G$ be a finitely generated $BCP$-group. Then $G$ has only a finite number of normal
subgroups of finite index.
\endproclaim

\demo{Proof} Suppose that $G$ is $m$-generated. First we claim that the order of each finite quotient of $G$ is  bounded by some fixed integer, say $f$.

Indeed, let $G/M$ be a finite quotient of $G$. Since $G$ is a $BCP$-group, it follows that $G$ is a $CP$-group and so are also the finite quotients $G/M$ of $G$.
By  Theorem 4 in [3], the order of each finite $CP$-group has a bounded number of prime divisors. Denote this bound by $d$.  Thus all finite quotients
$G/M$ of $G$ satisfy $|\pi(G/M)|\leq d$ and suppose that $|\pi(G/N)|$ is maximal among all finite quotients of $G$.
If some finite quotient $G/M$ of $G$ contains an element of prime order $p$
and  $p\notin \pi(G/N)$, then consider the quotient $G/S$, where $S=M\cap N$. Then $G/S$ is a finite quotient of $G$, such that $p\in \pi(G/S)$
and $\pi(G/N)\subset \pi(G/S)$, in contradiction to the maximality of $|\pi(G/N)|$. Hence, for each finite quotient
$G/M$ of $G$, the set $\pi(G/M)$ is a subset of $\pi(G/N)$.
Since $G/N$ is a  $BCP$-group, it follows that
$$\exp(G/N)\leq t=\prod_{p\in \pi(G/N)}p^{u_p},$$
and since $G/N$ is a finite group, $t$ is a finite integer.
Therefore
  $\exp(G/M)\leq t$ for all finite quotients $G/M$ of $G$.
Since each such finite quotient is $m$-generated and of exponent $\leq t$, it follows by the
Zelmanov positive solution of the Restricted Burnside Problem (see [11] and [12]) that their order is bounded by some fixed integer, say $f$, as claimed.

Since $G$ is finitely generated, there are only a finite number of normal subgroups $M$ of $G$ with a given finite index.
Since that index is bounded by $f$,
it follows that there exist only finitely many normal subgroups of $G$ of finite index.
\qed
\enddemo

Theorem 1 will be applied in the proofs of the next three theorem and indirectly also in the proof of Theorem 7.

\proclaim {Theorem 2} Let $G$ be a finitely generated residually finite $BCP$-group. Then $G$ is a finite group.
\endproclaim
\demo{Proof}  Since $G$ is residually
finite, for each non-trivial element $g\in G$  there exists a normal subgroup $M(g)$ of $G$ such that $g\notin M(g)$ and $G/M(g)$ is finite.
Let $T$ denote the intersection of the groups $M(g)$ for all non-trivial elements $g$ of $G$. Since $G$ is a finitely generated $BCP$-group, it follows by
Theorem 1 that there exist only finitely many normal subgroups of $G$ of finite index.  Therefore
$G/T$ is a finite group. But for each non-trivial $g\in G$ we have $g\notin M(g)$ , so $T=\{1\}$  and $G$ is a finite group, as required.
\qed
\enddemo

It is well known that the  residually finite property is inherited  by subgroups. This result follows from the fact that if $H$ and $M$ are
subgroups of a group $G$, then $[G:M]\geq [H:H\cap M]$. Using this fact and Theorem 2, we obtain the  following theorem.

\proclaim {Theorem 3} Let $G$ be a residually finite $BCP$-group. Then $G$ is a locally finite group.
\endproclaim
\demo {Proof} Let $H$ be a finitely generated subgroup of $G$. Then $H$ is a finitely generated residually finite $BCP$-group and hence it is
finite by Theorem 2. Thus $G$ is a locally finite group, as required.
\qed
\enddemo

\proclaim {Theorem 4} Let $G$ be a locally graded $BCP$-group. Then $G$ is a locally finite group. Therefore a finitely generated locally graded $BCP$-group is a finite group.
\endproclaim

\demo {Proof} Since $G$ is a locally graded group, each non-trivial finitely generated subgroup of $G$ has a proper normal subgroup of finite index.
Let $H$ be a finitely generated subgroup of $G$ and let $N$ be the intersection of all
normal subgroups of $H$ of finite index. Clearly $N$ is a normal subgroup of $H$. Since $H$ is a finitely generated $BCP$-group, it follows by Theorem 1 that $H$
has only a finite number of normal subgroups of finite index. Therefore $H/N$ is a finite group and
$N$ is a finitely generated subgroup of $G$. Since $G$ is locally graded, if $N$ is non-trivial, then it has a proper normal subgroup $T$ of finite index. Hence $T$ is also of finite index in $H$ and  it contains a subgroup $S$ normal in $H$ and of finite index in $H$. Thus we have $N \leq S\leq T < N$, a contradiction. So $N$ is trivial and $H$ is finite. Therefore $G$ is a locally finite group, as required.
\qed
\enddemo

Finally, we shall deal with $BCP$-groups satisfying the condition $2,3\in \pi(G)$. We shall prove the following result.

\proclaim {Theorem 5} Let $G$ be a  $BCP$-group satisfying $2,3\in \pi(G)$ and suppose that $u_2=1$ and $u_3\in \{1,2\}$. Then $G$ is a locally finite group.
\endproclaim

\demo{Proof}  Since $G$ is a periodic group and $2\in \pi(G)$, it follows that $G$ contains an involution. Let $t$ be any involution in $G$. Since $G$ is a $BCP$-group, $C_G(t)$
is a $2$-subgroup of $G$ and since $u_2=1$, it follows that $C_G(t)$ is an elementary abelian $2$-subgroup of $G$. As $G$ is a periodic group, Theorem 2(2) in V.D. Mazurov's paper
[6] implies that
one of the following statements holds:

(2.1) $G=A\langle t\rangle$, where $A$ is an abelian periodic subgroup of $G$ without involutions, and $a^t=a^{-1}$ for every $a\in A$.

(2.2) $G$ is an extension of an abelian $2$-group by a group without involutions.

(2.3) $G$ is isomorphic to $PGL_2(P)$, where $P$ is a locally finite field of characteristic $2$.

If (2.1) holds, then $A$ is an abelian periodic normal subgroup of $G$. Since $A$ and $G/A$ are locally finite, it follows  by the Schmidt's theorem (see 14.3.1 in
[9]) that $G$ is locally finite, as required.

If (2.3) holds, then $P$ being a locally finite field implies that $G$ is locally finite, as required.

It remains to deal with the case (2.2). In this case, there exists a normal elementary abelian $2$-subgroup $T$ of $G$, such that $G/T$ is a periodic
group with no involutions. Since $G$ is a $BCP$-group, it follows that $C_G(T)=T$ and hence $G/T$ is a periodic
subgroup of $Aut(T)$ without involutions. Moreover,  $o(gT)=o(g)$ for every non-trivial element $gT$ of $G/T$ and $G/T$ acts fixed point freely on $T$. By our assumptions
$G/T$ contains an element of order $3$ and by Lemma 1 in Zhurtov and Mazurov's  paper [13], each element of $G/T$ of order $3$ is in the center of $G/T$. Since $G/T$ is also a $BCP$-group,
it follows that $G/T$ is a $3$-group. If $u_3=1$, then $G/T$ is of exponent $3$ and hence it is abelian. Suppose finally that $u_3=2$ and $G/T$ is of exponent $9$.
Since every element of order $3$ in $G/T$ is in the center of $G/T$, it follows that  $(G/T)/(Z(G/T))$ is of exponent $3$ and by Lemmas 12.3.5 and 12.3.6 in [9], $(G/T)/(Z(G/T))$ is
a nilpotent group. Therefore $G/T$ is a periodic nilpotent group, and it follows by 5.2.18 in [9]
 that $G/T$ is locally finite. Since $T$ is a periodic abelian group, it
is also locally finite and by the Schmidt's theorem $G$ is locally finite, as required.

The proof of Theorem 5 is now complete.
\qed
\enddemo

\heading III. $BSP$-groups \\
\endheading

Finally, we shall deal with $BSP$-groups. Since each $BSP$-group is a $BCP$-group, all the results of Section II are valid for
$BSP$-groups as well.

The definition of the $BSP$-groups enables us to prove the following result, which does not hold for $BCP$-groups.

\proclaim {Theorem 6} Let $G$ be a locally finite $BSP$-group. Then $G$ is a finite group.
\endproclaim
\demo {Proof} Since $G$ is a locally finite $BSP$-group, it follows by the Main Theorem of [3] that $|\pi(G)|$ is bounded.
If $X$ is a finite subset of $G$, then
$$|\langle X\rangle|\leq \prod_{p\in \pi(G)}p^{v_p}.$$
Since $\prod_{p\in \pi(G)}p^{v_p}$ is a finite integer,
it follows that $G$ is a finite group, as required.
\qed
\enddemo

This theorem does not hold for $BCP$-group, since if $p$ is a prime, then an infinite abelian $p$-group of finite exponent is
a locally finite $BCP$-group.

The main result of this section is the following strengthening of Theorem 4 for $BSP$-groups.

\proclaim {Theorem 7} Let $G$ be a  locally graded $BSP$-group. Then $G$ is a finite group.
\endproclaim

\demo {Proof} By Theorem 4 applied to $BSP$-groups, $G$ is a locally finite group. Hence, by Theorem 6, $G$ is a finite group, as required.
\qed
\enddemo

Finally, we shall deal with $BSP$-groups satisfying the condition $2\in \pi(G)$.
First we prove the following theorem.

\proclaim {Theorem 8} Let $G$ be a $BSP$-$2$-group. Then $G$ is a finite group.
\endproclaim

\demo{Proof} If $G$ is an infinite $2$-group and $K$ is a finite subgroup of $G$, then by Theorem 14.4.1 in [9] $N_G(K)>K$.
If $G$ is also a $BSP$-group, then it is periodic, and it follows that there exists an infinite  series of finite $2$-subgroups of $G$
with increasing orders, in contradiction to the definition of a $BSP$-group. Hence a $BSP$-$2$-group is a finite group.
\qed
\enddemo

Our final main result is the following theorem.

\proclaim {Theorem 9} Let $G$ be a  $BSP$-group satisfying $2\in \pi(G)$. Then $G$ is a locally finite group.
\endproclaim

\demo{Proof} Let $t$ be an involution in $G$. Since $G$ is a $BSP$-group, $C_G(t)$ is a $BSP$-$2$-group and it is finite by Theorem 8.
Since that is true for any involution in $G$, it follows by Corollary 2 in the paper [10] of V.P. Shunkov that $G$ is locally finite, as claimed.
\qed
\enddemo

\enddocument